\documentclass[11pt]{article}
\usepackage{latexsym,bm}
\usepackage{color}
\usepackage{amssymb, graphicx, amsmath}  

\def \[{\begin{equation}}
\def \]{\end{equation}}

\newtheorem{theorem}{Theorem}[section]

\newtheorem{definition}{Definition}[section]

\newtheorem{lemma}{Lemma}[section]

\newtheorem{remark}{Remark}

\begin{document}

\begin{center}
{\bf The first order partial differential equations resolved with any derivatives}

\medskip

  {\bf  Jianfeng Wang }

 Department of Physics, Stern College at Yeshiva University,\\ 245 Lexington Avenue, New York, NY10016, USA.\\Department of Mathematics, Hohai University, Nanjing, 210098,
  P.R. China.

\end{center}
\bigskip
{\narrower \noindent  {\bf Abstract.} In this paper we discuss the first order partial differential equations resolved with any derivatives. At first, we transform the first order partial differential equation resolved with respect to a time derivative into a system of linear equations. Secondly, we convert it into a system of the first order linear partial
differential equations with constant coefficients and nonlinear algebraic
equations. Thirdly, we solve them by the Fourier transform and convert them into the equivalent integral equations. At last, we extend to discuss the first order partial differential equations resolved with respect to time derivatives and the general scenario resolved with any derivatives.\\
\noindent{\bf Keywords.} linear algebra, Fourier transform, integral equations.
    \par }

\vskip 1.0 true cm
 \section{Introduction}
The goal of this paper is to transform the first order partial differential equations resolved with any derivatives as follows into the integral equations,
\[v_{j}=f_{j}(v_{m+1},\ v_{m+2},\cdots,\ v_{5m},\ x,\ y,\ z,\ t),\ 1\leq j\leq m,\] where\ $t\in [0,\ T],\ (x,\ y,\ z)^{T}\in \Omega\subset R^{3}$, and\ $\Omega$\ is bounded,\ $\partial\Omega$\ is smooth or piecewise smooth, \[u=(u_{1}(x,\ y,\ z,\ t),\cdots,u_{m}(x,\ y,\ z,\ t))^{T}\in C^{1}(\Omega\times (0,\ T)),\] the initial conditions and boundary conditions are\[u|_{t=0}\in L^{2}(\Omega),\ u|_{\partial\Omega\times(0,\ T)}\in L^{2}(\partial\Omega\times(0,\ T)),\] $f_{j},\ 1\leq j\leq m,$\ are continuous functions,\ $v_{1},\ v_{2},\cdots,\ v_{5m}$\ is a permutation of the components of\ $u_{t},\ u,\ u_{x},\ u_{y},\ u_{z}.$\\After we did it, most of the Mathematical-Physics equation in [1] as follows can be transformed into the integral equations.
\\ 1. Eikonal equation\[ |Du|=1.\]2. Nonlinear Poisson eqation\[ -\triangle u=f(u).\]3. p-Laplacian equation\[ div(|Du|^{p-2}Du)=0.\]
4. Minimal surface equation\[ div\left(\cfrac{Du}{(1+|Du|^{2})^{1/2}}\right)=0.\] 5. Monge-Amp$\grave{e}$re equation\[ det(D^{2}u)=f.\]
6. Hamilton-Jacobi equation\[u_{t}+H(Du,\ x)=0.\]7. Scalar conservation law\[u_{t}+divF(u)=0.\]8. Inviscid Burgers' equation\[ u_{t}+uu_{x}=0.\]
9. Scalar reaction-diffusion equation\[u_{t}-\triangle u=f(u).\]10. Porous medium equation\[ u_{t}-\triangle(u^{\gamma})=0.\]
11. Nonlinear wave equations\[\bf{ u_{tt}-\triangle u=f(u),\ u_{tt}-div\ a(Du)=0.}\]12. Korteweg-de Vries (KdV) equation\[ u_{t}+uu_{x}+u_{xxx}=0.\]
13. System of conservation law\[\bf u_{t}+divF(u)=0.\]14. Reaction-diffusion system\[\bf u_{t}-\triangle u=f(u).\]
From\ $f_{j},\ 1\leq j\leq m,$\ are continuous functions, we can't use the Cauchy-Kovalevskaya theorem and the linear methods such as elliptic, parabola and hyperbola. In this paper, we will gradually transform Eq.(1.1) into the equivalent integral equations. In order to read easily, we discuss the first order partial differential equation resolved with respect to a time derivative at first.
 \section{One equation} \setcounter{equation}{0}
 In this section, we discuss the first order partial differential equation resolved with respect to a time derivative as follows,
\[u_{t}=f(u,\ u_{x},\ u_{y},\ u_{z},\ x,\ y,\ z,\ t), \]where\ $t\in [0,\ T],\ f$\ is a continuous function,\ $(x,\ y,\ z)^{T}\in \Omega\subset R^{3}$, and\ $\Omega$\ is bounded,\ $\partial\Omega$\ is smooth or piecewise smooth. For simplicity, we assume\ $u=u(x,\ y,\ z,\ t)\in C^{1}(\Omega\times (0,\ T))$.\ The initial condition and boundary condition are\[u|_{t=0}\in L^{2}(\Omega),\ u|_{\partial\Omega\times(0,\ T)}\in L^{2}(\partial\Omega\times(0,\ T)).\]At first we transform Eq.(2.1) into the
linear equations on unknown functions as we have done in [2],
 \[u_{t}-au-bu_{x}-cu_{y}-du_{z}-v=0,\]where\ $a,\ b,\ c,\ d$\ are real constants to be determined,
 \[v=f(u,\ u_{x},\ u_{y},\ u_{z},\ x,\ y,\ z,\ t)-au-bu_{x}-cu_{y}-du_{z}.\] Let's introduce\ $X=(x_1,\ x_2,\  x_3,\  x_4,\  x_5,\  x_6)^{T}$, where\[x_1=u_{t},\ x_2=u,\ x_3=u_{x},\ x_4=u_{y},\ x_5=u_{z},\ x_6=v.\]Then Eq.(2.1) is equivalent to\[\alpha^{T}X=0,\] where\ $\alpha=(1,\ -a,\ -b,\ -c,\ -d,\ -1)^{T}$.\\
 After we solve Eq.(2.6) by linear algebra, we obtain\ $\exists$\ independent variable vector\ $Z$, such that\[X=\left(
                                                                                                                                     \begin{array}{c}
                                                                                                                                       \beta^{T} \\
                                                                                                                                       E \\
                                                                                                                                     \end{array}
                                                                                                                                   \right)Z,\ \mbox{where}\ \beta=(a,\ b,\ c,\ d,\ 1)^{T}.\]
                                                                                                                                  We should discuss the independent variable vector\ $Z$\ as follows,
 \[u_{t}=\beta^{T}Z,\ u=e_{1}^{T}Z,\ u_{x}=e_{2}^{T}Z,\ u_{y}=e_{3}^{T}Z,\ u_{z}=e_{4}^{T}Z,\ v=e_{5}^{T}Z.\]
 And we obtain Eq.(2.1) is equivalent to the following system respect to\ $Z$,\begin{eqnarray}&&\cfrac{\partial e_{1}^{T}Z}{\partial t}=\beta^{T}Z,\\&& \cfrac{\partial e_{1}^{T}Z}{\partial x}=e_{2}^{T}Z,\\&& \cfrac{\partial e_{1}^{T}Z}{\partial y}=e_{3}^{T}Z,\\&& \cfrac{\partial e_{1}^{T}Z}{\partial z}=e_{4}^{T}Z,\\&& f(e_{1}^{T}Z,\ e_{2}^{T}Z,\ e_{3}^{T}Z,\ e_{4}^{T}Z,\ x,\ y,\ z,\ t)-ae_{1}^{T}Z-be_{2}^{T}Z-ce_{3}^{T}Z-de_{4}^{T}Z=e_{5}^{T}Z.\end{eqnarray}
     \begin{remark} \label{remark1}In fact, we should assume\[X=XI_{\Omega\times(0,\ T)}(x,\ y,\ z,\ t),\ \mbox{or}\ Z=ZI_{\Omega\times(0,\ T)}(x,\ y,\ z,\ t),\] if Eq.(2.1) is only satisfied on\ $\Omega\times(0,\ T)$.\end{remark}
 We will solve Eq.(2.9) to Eq.(2.12) by the Fourier transform as follows. Then Eq.(2.13) is the goal. \\ \begin{eqnarray}\int_{0}^{T}dt\int_{\Omega}e^{-it\xi_{0}-ix\xi_{1}-iy\xi_{2}-iz\xi_{3}}\cfrac{\partial e_{1}^{T}Z}{\partial t}dxdydz=\int_{0}^{T}dt\int_{\Omega}e^{-it\xi_{0}-ix\xi_{1}-iy\xi_{2}-iz\xi_{3}}\beta^{T}Zdxdydz,&&\\ \int_{0}^{T}dt\int_{\Omega}e^{-it\xi_{0}-ix\xi_{1}-iy\xi_{2}-iz\xi_{3}} \cfrac{\partial e_{1}^{T}Z}{\partial x}dxdydz=\int_{0}^{T}dt\int_{\Omega}e^{-it\xi_{0}-ix\xi_{1}-iy\xi_{2}-iz\xi_{3}}e_{2}^{T}Zdxdydz,&&\\ \int_{0}^{T}dt\int_{\Omega}e^{-it\xi_{0}-ix\xi_{1}-iy\xi_{2}-iz\xi_{3}}\cfrac{\partial e_{1}^{T}Z}{\partial y}dxdydz=\int_{0}^{T}dt\int_{\Omega}e^{-it\xi_{0}-ix\xi_{1}-iy\xi_{2}-iz\xi_{3}}e_{3}^{T}Zdxdydz,&&\\ \int_{0}^{T}dt\int_{\Omega}e^{-it\xi_{0}-ix\xi_{1}-iy\xi_{2}-iz\xi_{3}} \cfrac{\partial e_{1}^{T}Z}{\partial z}dxdydz=\int_{0}^{T}dt\int_{\Omega}e^{-it\xi_{0}-ix\xi_{1}-iy\xi_{2}-iz\xi_{3}}e_{4}^{T}Zdxdydz.\end{eqnarray}\\
To denote easily, we define the Fourier transform on\ $\Omega\times(0,\ T)$\ as
follows.\begin{definition}\label{definition} $\forall\ f(x,\ y,\ z,\ t)\in L^{2}(\Omega\times(0,\ T))$,
\begin{eqnarray}FI(f(x,\ y,\ z,\ t))&=&\int_{0}^{T}\int_{\Omega}f(x,\ y,\ z,\ t)e^{-it\xi_{0}-ix\xi_{1}-iy\xi_{2}-iz\xi_{3}}dxdydzdt\\&=&F(f(x,\ y,\ z,\ t)I_{\Omega\times(0,\ T)}(x,\ y,\ z,\ t)),\end{eqnarray}
where\ $F$\ means the Fourier transform and\ $I_{\Omega\times(0,\ T)}(x,\ y,\ z,\ t)$\
is the characteristic function. In the following, we write\ $I_{\Omega\times(0,\ T)}(x,\ y,\ z,\ t)$\ into\ $I_{\Omega\times(0,\ T)}$.\end{definition}Then we obtain\begin{eqnarray}
FI(\cfrac{\partial e_{1}^{T}Z}{\partial t})&=&\int_{0}^{T}dt\int_{\Omega}e^{-it\xi_{0}-ix\xi_{1}-iy\xi_{2}-iz\xi_{3}}(\cfrac{\partial e_{1}^{T}Z}{\partial t})dxdydz
\\&=&\int_{\Omega}(e_{1}^{T}Z)e^{-it\xi_{0}}|_{t=0}^{t=T}e^{-ix\xi_{1}-iy\xi_{2}-iz\xi_{3}}dxdydz
+\\&&i\xi_{0}\int_{0}^{T}dt\int_{\Omega}e^{-it\xi_{0}-ix\xi_{1}-iy\xi_{2}-iz\xi_{3}}(e_{1}^{T}Z)dxdydz
\\&=&f_{0}+ i\xi_{0} FI(e_{1}^{T}Z),\end{eqnarray}\begin{eqnarray}
FI(\cfrac{\partial e_{1}^{T}Z}{\partial x})&=&\int_{0}^{T}dt\int_{\Omega}e^{-it\xi_{0}-ix\xi_{1}-iy\xi_{2}-iz\xi_{3}}(\cfrac{\partial e_{1}^{T}Z}{\partial x})dxdydz
\\&=&\int_{0}^{T}\int_{\partial\Omega}(e_{1}^{T}Z)n_{1}e^{-it\xi_{0}-ix\xi_{1}-iy\xi_{2}-iz\xi_{3}}dxdydz
+\\&&i\xi_{1}\int_{0}^{T}dt\int_{\Omega}e^{-it\xi_{0}-ix\xi_{1}-iy\xi_{2}-iz\xi_{3}}(e_{1}^{T}Z)dxdydz
\\&=&f_{1}+ i\xi_{1} FI(e_{1}^{T}Z),\end{eqnarray} \begin{eqnarray}
FI(\cfrac{\partial e_{1}^{T}Z}{\partial y})&=&\int_{0}^{T}dt\int_{\Omega}e^{-it\xi_{0}-ix\xi_{1}-iy\xi_{2}-iz\xi_{3}}(\cfrac{\partial e_{1}^{T}Z}{\partial y})dxdydz
\\&=&\int_{0}^{T}\int_{\partial\Omega}(e_{1}^{T}Z)n_{2}e^{-it\xi_{0}-ix\xi_{1}-iy\xi_{2}-iz\xi_{3}}dxdydz
+\\&&i\xi_{2}\int_{0}^{T}dt\int_{\Omega}e^{-it\xi_{0}-ix\xi_{1}-iy\xi_{2}-iz\xi_{3}}(e_{1}^{T}Z)dxdydz
\\&=&f_{2}+ i\xi_{2} FI(e_{1}^{T}Z),\end{eqnarray} \begin{eqnarray}
FI(\cfrac{\partial e_{1}^{T}Z}{\partial z})&=&\int_{0}^{T}dt\int_{\Omega}e^{-it\xi_{0}-ix\xi_{1}-iy\xi_{2}-iz\xi_{3}}(\cfrac{\partial e_{1}^{T}Z}{\partial z})dxdydz
\\&=&\int_{0}^{T}\int_{\partial\Omega}(e_{1}^{T}Z)n_{3}e^{-it\xi_{0}-ix\xi_{1}-iy\xi_{2}-iz\xi_{3}}dxdydz
+\\&&i\xi_{3}\int_{0}^{T}dt\int_{\Omega}e^{-it\xi_{0}-ix\xi_{1}-iy\xi_{2}-iz\xi_{3}}(e_{1}^{T}Z)dxdydz
\\&=&f_{3}+ i\xi_{3} FI(e_{1}^{T}Z),\end{eqnarray}where\begin{eqnarray}f_{0}&=&\int_{\Omega}(A_{2}e^{-iT\xi_{0}}-A_{1})e^{-ix\xi_{1}-iy\xi_{2}-iz\xi_{3}}dxdydz,\\
f_{1}&=&\int_{0}^{T}dt\int_{\partial \Omega}A_{3}n_{1}e^{-it\xi_{0}-ix\xi_{1}-iy\xi_{2}-iz\xi_{3}}dS,\\
f_{2}&=&\int_{0}^{T}dt\int_{\partial \Omega}A_{3}n_{2}e^{-it\xi_{0}-ix\xi_{1}-iy\xi_{2}-iz\xi_{3}}dS,\\
f_{3}&=&\int_{0}^{T}dt\int_{\partial \Omega}A_{3}n_{3}e^{-it\xi_{0}-ix\xi_{1}-iy\xi_{2}-iz\xi_{3}}dS.\\
A_{1}&=&u|_{t=0},\ A_{2}=u|_{t=T},\ A_{3}=u|_{\partial\Omega\times(0,\ T)},
\end{eqnarray}\ $n_{k}$\ is the\ $k$th component of
the normal vector to\ $\partial \Omega,\ k=1,\ 2,\ 3$. We only need\[A_{1},\ A_{2}\in L^{2}(\Omega),\ A_{3}\in L^{2}(\partial\Omega\times(0,\ T)).\]
Now we transformed Eq.(2.9) to Eq.(2.12) into the following.\[\label
{A-Problem}BFI(Z)=\beta_{1},\]where\[B=\left(
                                         \begin{array}{c}
                                        i\xi_{0}e_{1}^{T}-\beta^{T} \\
                                         i\xi_{1}e_{1}^{T}-e_{2}^{T} \\
                                         i\xi_{2}e_{1}^{T}-e_{3}^{T} \\
                                         i\xi_{3}e_{1}^{T}-e_{4}^{T} \\
                                         \end{array}
                                       \right)_{4\times5}
=(B_{1},\ -B_{2}),\]
\[B_{1}=\left(
          \begin{array}{cccc}
            i\xi_{0}-a, & -b, &-c, &-d\\
            i\xi_{1}, & -1, & 0, & 0 \\
            i\xi_{2}, & 0, &-1, & 0 \\
            i\xi_{3}, & 0, &0, & -1 \\
          \end{array}
        \right)_{4\times4},\] \[B_{2}=(1,\ 0,\ 0,\ 0)^{T},\ \beta_{1}=(-f_{0},\ -f_{1},\ -f_{2},\ -f_{3})^{T}.\]
If we assume
\[ Z=\left(
       \begin{array}{c}
         Z_{1} \\
         Z_{2} \\
       \end{array}
     \right),\ \mbox{where}\ Z_{1}\ \mbox{is the first\ $4$\ componenets of}\ Z,\] then we obtain\[B_{1}FI(Z_{1})=\beta_{1}+B_{2}FI(Z_{2}).\]It is not very difficult to work out\[ det(B_{1})=a+b i\xi_{1}+c i\xi_{2}+d i\xi_{3}- i\xi_{0}=a_{1},\] \[B_{1}^{-1}=-\left(
                                                                                                           \begin{array}{cccc}
                                                                                                             a_{1}^{-1}, & a_{1}^{-1}b, &a_{1}^{-1}c, & a_{1}^{-1}d \\
                                                                                                            i\xi_{1}a_{1}^{-1}, &  1+i\xi_{1}a_{1}^{-1}b, & i\xi_{1}a_{1}^{-1}c, &  i\xi_{1}a_{1}^{-1}d \\
                                                                                                              i\xi_{2}a_{1}^{-1}, &  i\xi_{2}a_{1}^{-1}b, & 1+i\xi_{2}a_{1}^{-1}c, &  i\xi_{2}a_{1}^{-1}d \\
                                                                                                             i\xi_{3}a_{1}^{-1}, &  i\xi_{3}a_{1}^{-1}b, & i\xi_{3}a_{1}^{-1}c, &  1+i\xi_{3}a_{1}^{-1}d \\
                                                                                                           \end{array}
                                                                                                         \right),\ B_{1}^{-1}B_{2}=-\left(
                                                                                                                                     \begin{array}{c}
                                                                                                                                      a_{1}^{-1} \\
                                                                                                                                       i\xi_{1}a_{1}^{-1} \\
                                                                                                                                       i\xi_{2}a_{1}^{-1} \\
                                                                                                                                       i\xi_{3}a_{1}^{-1} \\
                                                                                                                                     \end{array}
                                                                                                                                   \right).\]
                                                                                                                                   If we assume\ $C=\{\xi^{\prime}|a_{1}=0\}$, where\ $\xi^{\prime}=(\xi_{0},\ \xi_{1},\ \xi_{2},\ \xi_{3})^{T}$, then the measure of\ $C$\ is\ $0$. And we obtain
 \[FI(Z_{1})(1-I_{C}(\xi^{\prime}))=B_{1}^{-1}\beta_{1}(1-I_{C}(\xi^{\prime}))+B_{1}^{-1}B_{2}FI(Z_{2})(1-I_{C}(\xi^{\prime})).\]
We need some lemmas.
\begin{lemma} \label{lemma1}(Plancherel Theorem) If\ $f(x,\ y,\ z,\ t)\in L^{2}(R^{4})$, then\ $F(f(x,\ y,\ z,\ t))$\ exists, moreover\\(1)$\parallel F(f(x,\ y,\ z,\ t))\parallel_{L^{2}}=\parallel f(x,\ y,\ z,\ t)\parallel_{L^{2}}$,\\(2)$F^{-1}[F(f(x,\ y,\ z,\ t))]=f(x,\ y,\ z,\ t)$.\end{lemma}
\begin{lemma} \label{lemma1} If\ $f(x,\ y,\ z,\ t)\in L^{2}(R^{4})$,\ $C\subset R^{4}$, the measure of\ $C$\ is\ $0$, then \[F^{-1}([F(f(x,\ y,\ z,\ t))](1-I_{C}(\xi^{\prime})))=f(x,\ y,\ z,\ t).\]\end{lemma}
{\it Proof of lemma 2.2}. From the lemma 2.1, we know\ $F(f(x,\ y,\ z,\ t))\in L^{2}(R^{4})$. Therefore,
\[\int_{C}F(f(x,\ y,\ z,\ t))e^{it\xi_{0}+ix\xi_{1}+iy\xi_{2}+iz\xi_{3}}d\xi_{0}d\xi_{1}d\xi_{2}d\xi_{3}=0.\]And we obtain
\[ F^{-1}([F(f(x,\ y,\ z,\ t))](1-I_{C}(\xi^{\prime})))=F^{-1}[F(f(x,\ y,\ z,\ t))]=f(x,\ y,\ z,\ t).\]\\
From these two lemmas, we obtain\[F^{-1}[FI(Z_{1})(1-I_{C}(\xi^{\prime}))]=Z_{1}I_{\Omega\times(0,\ T)},\ F^{-1}[FI(Z_{2})(1-I_{C}(\xi^{\prime}))]=Z_{2}I_{\Omega\times(0,\ T)}.\]
Now we determine the parameters. We choose\ $a,\ b,\ c,\ d$, such that\ $F^{-1}(a_{1}^{-1})$\ exists. Then\ $F^{-1}(B_{1}^{-1}B_{2})$\ exists. And\ $F^{-1}[B_{1}^{-1}\beta_{1}(1-I_{C}(\xi^{\prime}))]\ \mbox{exists}.$\\
If we assume\[w_{1}(x,\ y,\ z,\ t)=F^{-1}[B_{1}^{-1}\beta_{1}(1-I_{C}(\xi^{\prime}))],\ w_{2}(x,\ y,\ z,\ t)=F^{-1}(B_{1}^{-1}B_{2}),\]
then we obtain\[Z_{1}I_{\Omega\times(0,\ T)}=w_{1}(x,\ y,\ z,\ t)+w_{2}(x,\ y,\ z,\ t)\ast(Z_{2}I_{\Omega\times(0,\ T)}),\]
where\begin{eqnarray}Z_{1}=(e_{1}^{T}Z,\ e_{2}^{T}Z,\ e_{3}^{T}Z,\ e_{4}^{T}Z)^{T}=(u,\ u_{x},\ u_{y},\ u_{z})^{T},&&\end{eqnarray}\[Z_{2}=e_{5}^{T}Z=f(e_{1}^{T}Z,\ e_{2}^{T}Z,\ e_{3}^{T}Z,\ e_{4}^{T}Z,\ x,\ y,\ z,\ t)-ae_{1}^{T}Z-be_{2}^{T}Z-ce_{3}^{T}Z-de_{4}^{T}Z=v.\]It is obvious\ $\exists\ \psi$, such that\ $Z_{2}=\psi(Z_{1})$. Therefore, we obtain
\[ Z_{1}I_{\Omega\times(0,\ T)}=w_{1}(x,\ y,\ z,\ t)+w_{2}(x,\ y,\ z,\ t)\ast(\psi(Z_{1})I_{\Omega\times(0,\ T)}).\]
If\ $Z_{1}$\ satisfied Eq.(2.60), then we let\ $Z_{2}=\psi(Z_{1})$. We obtain\ $ B FI(Z)=\beta_{1},\ \alpha^{T}X=0,$\ on\ $\Omega\times(0,\ T)$. Therefore,\ $e_{1}^{T}Z$\ is the solution of Eq.(2.1)\ on\ $\Omega\times(0,\ T)$. Hence we arrive at \begin{theorem} \label{Theorem2-1}\ $ w_1,\ w_2,\ \psi,$\ as we described, then Eq.(2.1) is equivalent to Eq.(2.60). \end{theorem}
Maybe you will say we should know\ $A_{2}=u|_{t=T}$. In fact, we shouldn't. If we choose the parameter\ $a<0$, then in\ $F^{-1}(a_{1}^{-1}f_{0})$, we obtain\[\int_{-\infty}^{+\infty}\cfrac{e^{-iT\xi_{0}}a_{3}}{i\xi_{0}+a_{2}}\ e^{it\xi_{0}}d\xi_{0}=e^{-a_{2}(t-T)}a_{3}I_{\{t\geq T\}},\]
where\[a_{2}=-(a+b i\xi_{1}+c i\xi_{2}+d i\xi_{3}),\ a_{3}=\int_{\Omega}A_{2}e^{-ix\xi_{1}-iy\xi_{2}-iz\xi_{3}}dxdydz.\] $I_{\{t\geq T\}}$\ means\ $A_{2}$\ doesn't work on\ $\Omega\times(0,\ T)$. If we only discuss\ $u$\ on\ $\Omega\times(0,\ T)$, then we only need to know\ $A_{1}=u|_{t=0},\ A_{3}=u|_{\partial\Omega\times(0,\ T)}$\ in Eq.(2.60). Hence we take\ $f_{0}$\ as\[f_{0}=\int_{\Omega}(-A_{1})e^{-ix\xi_{1}-iy\xi_{2}-iz\xi_{3}}dxdydz.\]\\
Maybe you will also say Eq.(2.60) is the kind of Hammerstain. We take it at beginning, but we change the mind after we obtain\ $w_{2}(x,\ y,\ z,\ t)$\ with\ $I_{\{t\geq 0\}}$.
 \section{The equations} \setcounter{equation}{0}
 In this section, we discuss the first order partial differential equations resolved with respect to time derivatives as follows,
\[u_{t}=f(u,\ u_{x},\ u_{y},\ u_{z},\ x,\ y,\ z,\ t), \]where\ $t\in [0,\ T],\ f=(f_{1},\cdots,f_{m})^{T}$\ is continuous,\ $(x,\ y,\ z)^{T}\in \Omega\subset R^{3}$, and\ $\Omega$\ is bounded,\ $\partial\Omega$\ is smooth or piecewise smooth. For simplicity, we assume \[u=(u_{1}(x,\ y,\ z,\ t),\cdots,u_{m}(x,\ y,\ z,\ t))^{T}\in C^{1}(\Omega\times (0,\ T)).\] The initial conditions and boundary conditions are\[u|_{t=0}\in L^{2}(\Omega),\ u|_{\partial\Omega\times(0,\ T)}\in L^{2}(\partial\Omega\times(0,\ T)).\]
We transform Eq.(3.1) into the
linear equations on unknown functions as follows,
 \[u_{t}-Au-Bu_{x}-Cu_{y}-Du_{z}-v=0,\]where\ $A,\ B,\ C,\ D$\ are real constants\ $m\times m$\ matrices to be determined,
 \[v=f(u,\ u_{x},\ u_{y},\ u_{z},\ x,\ y,\ z,\ t)-Au-Bu_{x}-Cu_{y}-Du_{z}.\] Let's introduce\ $X=(x_1^{T},\ x_2^{T},\  x_3^{T},\  x_4^{T},\  x_5^{T},\  x_6^{T})^{T}$, where\[x_1=u_{t},\ x_2=u,\ x_3=u_{x},\ x_4=u_{y},\ x_5=u_{z},\ x_6=v.\]Then Eq.(3.14) is equivalent to\[\alpha^{T}X=0,\] where\ $\alpha^{T}=(E,\ -A,\ -B,\ -C,\ -D,\ -E)$.\\
  After we solve Eq.(3.7) by linear algebra, we obtain\ $\exists$\ independent variable vector\ $Z$, such that\[X=\left(
                                                                                                                                     \begin{array}{c}
                                                                                                                                       \beta^{T} \\
                                                                                                                                       E_{5m} \\
                                                                                                                                     \end{array}
                                                                                                                                   \right)Z,\ \mbox{where}\ \beta^{T}=(A,\ B,\ C,\ D,\ E).\]
                                                                                                                                  We should discuss the independent variable vector\ $Z$\ as follows,
 \begin{eqnarray}&&u_{t}=\beta^{T}Z,\ u=E_{1}^{T}Z,\ u_{x}=E_{2}^{T}Z,\ u_{y}=E_{3}^{T}Z,\ u_{z}=E_{4}^{T}Z,\ v=E_{5}^{T}Z,\\&& \mbox{where}\ E_{j}=(e_{m(j-1)+1},\ e_{m(j-1)+2},\cdots,\ e_{mj}),\ 1\leq j\leq 5.\end{eqnarray}
 And we obtain Eq.(3.4) is equivalent to the following system respect to\ $Z$,\begin{eqnarray}&&\cfrac{\partial E_{1}^{T}Z}{\partial t}=\beta^{T}Z,\\&& \cfrac{\partial E_{1}^{T}Z}{\partial x}=E_{2}^{T}Z,\\&& \cfrac{\partial E_{1}^{T}Z}{\partial y}=E_{3}^{T}Z,\\&& \cfrac{\partial E_{1}^{T}Z}{\partial z}=E_{4}^{T}Z,\end{eqnarray}\begin{eqnarray} f(E_{1}^{T}Z,\ E_{2}^{T}Z,\ E_{3}^{T}Z,\ E_{4}^{T}Z,\ x,\ y,\ z,\ t)-AE_{1}^{T}Z-BE_{2}^{T}Z-CE_{3}^{T}Z-DE_{4}^{T}Z=E_{5}^{T}Z.&&\end{eqnarray}
  We will solve Eq.(3.11) to Eq.(3.14) by the Fourier transform as follows. Then Eq.(3.15) is the goal. \\ \begin{eqnarray}\int_{0}^{T}dt\int_{\Omega}e^{-it\xi_{0}-ix\xi_{1}-iy\xi_{2}-iz\xi_{3}}\cfrac{\partial E_{1}^{T}Z}{\partial t}dxdydz=\int_{0}^{T}dt\int_{\Omega}e^{-it\xi_{0}-ix\xi_{1}-iy\xi_{2}-iz\xi_{3}}\beta^{T}Zdxdydz,&& \\ \int_{0}^{T}dt\int_{\Omega}e^{-it\xi_{0}-ix\xi_{1}-iy\xi_{2}-iz\xi_{3}} \cfrac{\partial E_{1}^{T}Z}{\partial x}dxdydz=\int_{0}^{T}dt\int_{\Omega}e^{-it\xi_{0}-ix\xi_{1}-iy\xi_{2}-iz\xi_{3}}E_{2}^{T}Zdxdydz,&&\\ \int_{0}^{T}dt\int_{\Omega}e^{-it\xi_{0}-ix\xi_{1}-iy\xi_{2}-iz\xi_{3}}\cfrac{\partial E_{1}^{T}Z}{\partial y}dxdydz=\int_{0}^{T}dt\int_{\Omega}e^{-it\xi_{0}-ix\xi_{1}-iy\xi_{2}-iz\xi_{3}}E_{3}^{T}Zdxdydz,&&\\ \int_{0}^{T}dt\int_{\Omega}e^{-it\xi_{0}-ix\xi_{1}-iy\xi_{2}-iz\xi_{3}} \cfrac{\partial E_{1}^{T}Z}{\partial z}dxdydz=\int_{0}^{T}dt\int_{\Omega}e^{-it\xi_{0}-ix\xi_{1}-iy\xi_{2}-iz\xi_{3}}E_{4}^{T}Zdxdydz.&&\end{eqnarray}\\
By using the definition 2.1, we obtain\begin{eqnarray}
FI(\cfrac{\partial E_{1}^{T}Z}{\partial t})&=&\int_{0}^{T}dt\int_{\Omega}e^{-it\xi_{0}-ix\xi_{1}-iy\xi_{2}-iz\xi_{3}}(\cfrac{\partial E_{1}^{T}Z}{\partial t})dxdydz
\\&=&\int_{\Omega}(E_{1}^{T}Z)e^{-it\xi_{0}}|_{t=0}^{t=T}e^{-ix\xi_{1}-iy\xi_{2}-iz\xi_{3}}dxdydz
+\\&&i\xi_{0}\int_{0}^{T}dt\int_{\Omega}e^{-it\xi_{0}-ix\xi_{1}-iy\xi_{2}-iz\xi_{3}}(E_{1}^{T}Z)dxdydz
\\&=&f_{0}+ i\xi_{0} FI(E_{1}^{T}Z),\end{eqnarray}\begin{eqnarray}
FI(\cfrac{\partial E_{1}^{T}Z}{\partial x})&=&\int_{0}^{T}dt\int_{\Omega}e^{-it\xi_{0}-ix\xi_{1}-iy\xi_{2}-iz\xi_{3}}(\cfrac{\partial E_{1}^{T}Z}{\partial x})dxdydz
\\&=&\int_{0}^{T}\int_{\partial\Omega}(E_{1}^{T}Z)n_{1}e^{-it\xi_{0}-ix\xi_{1}-iy\xi_{2}-iz\xi_{3}}dxdydz
+\\&&i\xi_{1}\int_{0}^{T}dt\int_{\Omega}e^{-it\xi_{0}-ix\xi_{1}-iy\xi_{2}-iz\xi_{3}}(E_{1}^{T}Z)dxdydz
\\&=&f_{1}+ i\xi_{1} FI(E_{1}^{T}Z),\end{eqnarray} \begin{eqnarray}
FI(\cfrac{\partial E_{1}^{T}Z}{\partial y})&=&\int_{0}^{T}dt\int_{\Omega}e^{-it\xi_{0}-ix\xi_{1}-iy\xi_{2}-iz\xi_{3}}(\cfrac{\partial E_{1}^{T}Z}{\partial y})dxdydz
\\&=&\int_{0}^{T}\int_{\partial\Omega}(E_{1}^{T}Z)n_{2}e^{-it\xi_{0}-ix\xi_{1}-iy\xi_{2}-iz\xi_{3}}dxdydz
+\\&&i\xi_{2}\int_{0}^{T}dt\int_{\Omega}e^{-it\xi_{0}-ix\xi_{1}-iy\xi_{2}-iz\xi_{3}}(E_{1}^{T}Z)dxdydz
\\&=&f_{2}+ i\xi_{2} FI(E_{1}^{T}Z),\end{eqnarray} \begin{eqnarray}
FI(\cfrac{\partial E_{1}^{T}Z}{\partial z})&=&\int_{0}^{T}dt\int_{\Omega}e^{-it\xi_{0}-ix\xi_{1}-iy\xi_{2}-iz\xi_{3}}(\cfrac{\partial E_{1}^{T}Z}{\partial z})dxdydz
\\&=&\int_{0}^{T}\int_{\partial\Omega}(E_{1}^{T}Z)n_{3}e^{-it\xi_{0}-ix\xi_{1}-iy\xi_{2}-iz\xi_{3}}dxdydz
+\\&&i\xi_{3}\int_{0}^{T}dt\int_{\Omega}e^{-it\xi_{0}-ix\xi_{1}-iy\xi_{2}-iz\xi_{3}}(E_{1}^{T}Z)dxdydz
\\&=&f_{3}+ i\xi_{3} FI(E_{1}^{T}Z),\end{eqnarray}where\begin{eqnarray}f_{0}&=&\int_{\Omega}(A_{2}e^{-iT\xi_{0}}-A_{1})e^{-ix\xi_{1}-iy\xi_{2}-iz\xi_{3}}dxdydz,\\
f_{1}&=&\int_{0}^{T}dt\int_{\partial \Omega}A_{3}n_{1}e^{-it\xi_{0}-ix\xi_{1}-iy\xi_{2}-iz\xi_{3}}dS,\\
f_{2}&=&\int_{0}^{T}dt\int_{\partial \Omega}A_{3}n_{2}e^{-it\xi_{0}-ix\xi_{1}-iy\xi_{2}-iz\xi_{3}}dS,\\
f_{3}&=&\int_{0}^{T}dt\int_{\partial \Omega}A_{3}n_{3}e^{-it\xi_{0}-ix\xi_{1}-iy\xi_{2}-iz\xi_{3}}dS.\\
A_{1}&=&u|_{t=0},\ A_{2}=u|_{t=T},\ A_{3}=u|_{\partial\Omega\times(0,\ T)},
\end{eqnarray}\ $n_{k}$\ is the\ $k$th component of
the normal vector to\ $\partial \Omega,\ k=1,\ 2,\ 3$. We only need\[A_{1},\ A_{2}\in L^{2}(\Omega),\ A_{3}\in L^{2}(\partial\Omega\times(0,\ T)).\]
Now we transformed the equations Eq.(3.11) to Eq.(3.14) into the following.\[\label
{A-Problem}BFI(Z)=\beta_{1},\]where\[B=\left(
                                         \begin{array}{c}
                                        i\xi_{0}E_{1}^{T}-\beta^{T} \\
                                         i\xi_{1}E_{1}^{T}-E_{2}^{T} \\
                                         i\xi_{2}E_{1}^{T}-E_{3}^{T} \\
                                         i\xi_{3}E_{1}^{T}-E_{4}^{T} \\
                                         \end{array}
                                       \right)_{4m\times5m}
=(B_{1},\ -B_{2}),\]
\[B_{1}=\left(
          \begin{array}{cccc}
            i\xi_{0}E-A, & -B, &-C, &-D\\
            i\xi_{1}E, & -E, & 0, & 0 \\
            i\xi_{2}E, & 0, &-E, & 0 \\
            i\xi_{3}E, & 0, &0, & -E \\
          \end{array}
        \right)_{4m\times4m},\] \[B_{2}=(E,\ 0,\ 0,\ 0)^{T},\ \beta_{1}=(-f_{0}^{T},\ -f_{1}^{T},\ -f_{2}^{T},\ -f_{3}^{T})^{T}.\]
If we assume
\[ Z=\left(
       \begin{array}{c}
         Z_{1} \\
         Z_{2} \\
       \end{array}
     \right),\ \mbox{where}\ Z_{1}\ \mbox{is the first\ $4m$\ componenets of}\ Z,\] then we obtain\[B_{1}FI(Z_{1})=\beta_{1}+B_{2}FI(Z_{2}).\]It is not very difficult to work out\[ det(B_{1})=det(A+ i\xi_{1}B+ i\xi_{2}C+ i\xi_{3}D- i\xi_{0}E)=det(B_{0}),\] where\[B_{0}=A+ i\xi_{1}B+ i\xi_{2}C+ i\xi_{3}D- i\xi_{0}E.\] \[B_{1}^{-1}=-\left(
                                                                                                           \begin{array}{cccc}
                                                                                                             B_{0}^{-1}, & B_{0}^{-1}B, &B_{0}^{-1}C, & B_{0}^{-1}D \\
                                                                                                            i\xi_{1}B_{0}^{-1}, &  E+i\xi_{1}B_{0}^{-1}B, & i\xi_{1}B_{0}^{-1}C, &  i\xi_{1}B_{0}^{-1}D \\
                                                                                                              i\xi_{2}B_{0}^{-1}, &  i\xi_{2}B_{0}^{-1}B, & E+i\xi_{2}B_{0}^{-1}C, &  i\xi_{2}B_{0}^{-1}D \\
                                                                                                             i\xi_{3}B_{0}^{-1}, &  i\xi_{3}B_{0}^{-1}B, & i\xi_{3}B_{0}^{-1}C, &  E+i\xi_{3}B_{0}^{-1}D \\
                                                                                                           \end{array}
                                                                                                         \right),\ B_{1}^{-1}B_{2}=-\left(
                                                                                                                                     \begin{array}{c}
                                                                                                                                      B_{0}^{-1} \\
                                                                                                                                       i\xi_{1}B_{0}^{-1} \\
                                                                                                                                       i\xi_{2}B_{0}^{-1} \\
                                                                                                                                       i\xi_{3}B_{0}^{-1} \\
                                                                                                                                     \end{array}
                                                                                                                                   \right).\]
                                                                                                                                   If we assume\ $C_{0}=\{\xi^{\prime}|det(B_{0})=0\}$, where\ $\xi^{\prime}=(\xi_{0},\ \xi_{1},\ \xi_{2},\ \xi_{3})^{T}$, then the measure of\ $C_{0}$\ is\ $0$. And we obtain
 \[FI(Z_{1})(1-I_{C_{0}}(\xi^{\prime}))=B_{1}^{-1}\beta_{1}(1-I_{C_{0}}(\xi^{\prime}))+B_{1}^{-1}B_{2}FI(Z_{2})(1-I_{C_{0}}(\xi^{\prime})).\]
From two lemmas in last section, we also obtain\[F^{-1}[FI(Z_{1})(1-I_{C_{0}}(\xi^{\prime}))]=Z_{1}I_{\Omega\times(0,\ T)},\ F^{-1}[FI(Z_{2})(1-I_{C_{0}}(\xi^{\prime}))]=Z_{2}I_{\Omega\times(0,\ T)}.\]
Now we determine the parameters. We choose\ $A,\ B,\ C,\ D$, such that\ $F^{-1}(B_{0}^{-1})$\ exists. Then\ $F^{-1}(B_{1}^{-1}B_{2})$\ exists. And\ $F^{-1}[B_{1}^{-1}\beta_{1}(1-I_{C_{0}}(\xi^{\prime}))]\ \mbox{exists}.$\\
If we assume\[w_{1}(x,\ y,\ z,\ t)=F^{-1}[B_{1}^{-1}\beta_{1}(1-I_{C_{0}}(\xi^{\prime}))],\ w_{2}(x,\ y,\ z,\ t)=F^{-1}(B_{1}^{-1}B_{2}),\]
then we obtain\[Z_{1}I_{\Omega\times(0,\ T)}=w_{1}(x,\ y,\ z,\ t)+w_{2}(x,\ y,\ z,\ t)\ast(Z_{2}I_{\Omega\times(0,\ T)}),\]
where\begin{eqnarray}Z_{1}&=&(E_{1}^{T}Z,\ E_{2}^{T}Z,\ E_{3}^{T}Z,\ E_{4}^{T}Z)^{T}=(u^{T},\ u_{x}^{T},\ u_{y}^{T},\ u_{z}^{T})^{T},\\Z_{2}&=&E_{5}^{T}Z=f(E_{1}^{T}Z,\ E_{2}^{T}Z,\ E_{3}^{T}Z,\ E_{4}^{T}Z,\ x,\ y,\ z,\ t)-\\&& AE_{1}^{T}Z-BE_{2}^{T}Z-CE_{3}^{T}Z-DE_{4}^{T}Z=v.\end{eqnarray}It is obvious\ $\exists\ \psi$, such that\ $Z_{2}=\psi(Z_{1})$. Therefore, we obtain
\[ Z_{1}I_{\Omega\times(0,\ T)}=w_{1}(x,\ y,\ z,\ t)+w_{2}(x,\ y,\ z,\ t)\ast(\psi(Z_{1})I_{\Omega\times(0,\ T)}).\]
If\ $Z_{1}$\ satisfied Eq.(3.58), then we let\ $Z_{2}=\psi(Z_{1})$. We obtain\ $ B FI(Z)=\beta_{1},\ \alpha^{T}X=0,$\ on\ $\Omega\times(0,\ T)$. Therefore,\ $E_{1}^{T}Z$\ is the solution of Eq.(3.1)\ on\ $\Omega\times(0,\ T)$. Hence we arrive at \begin{theorem} \label{Theorem2-1}\ $ w_1,\ w_2,\ \psi,$\ as we described, then Eq.(3.1) is equivalent to Eq.(3.58). \end{theorem}
Maybe you will say this time we should know\ $A_{2}=u|_{t=T}$. In fact, we also shouldn't. If we choose the parameter matrix\ $A=aE,\ a<0,\ B=bE,\ C=cE,\ D=dE,\ b,\ c,\ d\in R$, then in\ $F^{-1}(B_{0}^{-1}f_{0})$, we obtain\[\int_{-\infty}^{+\infty}(i\xi_{0}E+a_{2})^{-1}e^{-iT\xi_{0}}a_{3}\ e^{it\xi_{0}}d\xi_{0}=e^{-a_{2}(t-T)}a_{3}I_{\{t\geq T\}},\]
where\[ a_{2}=-(A+ i\xi_{1}B+ i\xi_{2}C+ i\xi_{3}D),\ a_{3}=\int_{\Omega}A_{2}e^{-ix\xi_{1}-iy\xi_{2}-iz\xi_{3}}dxdydz.\] $I_{\{t\geq T\}}$\ means\ $A_{2}$\ doesn't work on\ $\Omega\times(0,\ T)$. If we only discuss\ $u$\ on\ $\Omega\times(0,\ T)$, then we only need to know\ $A_{1}=u|_{t=0},\ A_{3}=u|_{\partial\Omega\times(0,\ T)}$\ in Eq.(3.58). Hence we take\ $f_{0}$\ as\[f_{0}=\int_{\Omega}(-A_{1})e^{-ix\xi_{1}-iy\xi_{2}-iz\xi_{3}}dxdydz.\]
\setcounter{remark}{0}\begin{remark} \label{remark1}Actually, the choice of the parameter matrix\ $A=aE,\ a<0,\ B=bE,\ C=cE,\ D=dE,\ b,\ c,\ d\in R$, is not unique. We only need to choose\ $A,\ B,\ C,\ D$\ which satisfy\ $Re(\lambda)>0$, where\ $\lambda$\ is the characteristic value of matrix\ $a_{2}$. We assume the following,\[a_{2}=PJP^{-1},\ J=\left(
                                 \begin{array}{ccc}
                                   J_{1} &  &  \\
                                    & \ddots &  \\
                                    &  &J_{\sigma} \\
                                 \end{array}
                               \right),\ J_{k}=\left(
                                                        \begin{array}{cccc}
                                                          \lambda_{k} & 1 &  &  \\
                                                          & \ddots & \ddots &  \\
                                                           &  & \ddots & 1 \\
                                                           &  &  & \lambda_{k} \\
                                                        \end{array}
                                                      \right)_{m_{k}\times m_{k}},\ 1\leq k\leq \sigma,\] where\ $P$\ is not related to\ $i\xi_{0},\ Re(\lambda_{k})>0,\ 1\leq k\leq \sigma.$\ Then we obtain
                                                      \begin{eqnarray}&&(i\xi_{0}E+a_{2})^{-1}=P(i\xi_{0}E+J)^{-1}P^{-1},\\&& (i\xi_{0}E+J)^{-1}=\left(
                                 \begin{array}{ccc}
                                   (i\xi_{0}E_{1}+J_{1})^{-1} &  &  \\
                                    & \ddots &  \\
                                    &  &(i\xi_{0}E_{\sigma}+J_{\sigma})^{-1} \\
                                 \end{array}
                               \right),\\&& (i\xi_{0}E_{k}+J_{k})^{-1}=\\&&\left(
                                                                      \begin{array}{ccccc}
                                                                       (i\xi_{0}+\lambda_{k})^{-1}, & -(i\xi_{0}+\lambda_{k})^{-2}, & (i\xi_{0}+\lambda_{k})^{-3}, & \cdots,& (-1)^{m_{k}-1}(i\xi_{0}+\lambda_{k})^{-m_{k}}, \\
                                                                        & \ddots & \ddots & \ddots & \vdots \\
                                                                         &  & \ddots &\ddots & (i\xi_{0}+\lambda_{k})^{-3} \\
                                                                         &  &  & \ddots & -(i\xi_{0}+\lambda_{k})^{-2} \\
                                                                         &  &  &  &(i\xi_{0}+\lambda_{k})^{-1}\\
                                                                      \end{array}
                                                                    \right)
                               _{m_{k}\times m_{k}}\end{eqnarray}where\ $1\leq k\leq \sigma.$\ Because\ $Re(\lambda_{k})>0,\ 1\leq k\leq \sigma$, we obtain
                               \[\int_{-\infty}^{+\infty}(i\xi_{0}+\lambda_{k})^{-1}e^{it\xi_{0}}d\xi_{0}=e^{-\lambda_{k}t}I_{\{t\geq 0\}},\ \int_{-\infty}^{+\infty}(i\xi_{0}+\lambda_{k})^{-n}e^{it\xi_{0}}d\xi_{0}=\cfrac{e^{-\lambda_{k}t}t^{n-1}}{(n-1)!}I_{\{t\geq 0\}},\ n\geq 1.\] Therefore, we obtain
                               \[\int_{-\infty}^{+\infty}(i\xi_{0}E+a_{2})^{-1}e^{it\xi_{0}}d\xi_{0}=e^{-a_{2}t}I_{\{t\geq 0\}}.\] Hence (3.59) stands.
\end{remark}
\section{General scenario}\setcounter{equation}{0}
In this section, we will consider the general first order partial differential equations. But what the general scenario should be? Is it the scenario as follows,\[f_{j}(u_{t},\ u,\ u_{x},\ u_{y},\ u_{z},\ x,\ y,\ z,\ t)=0,\ 1\leq j\leq m,\] where\ $t\in [0,\ T],\ (x,\ y,\ z)^{T}\in \Omega\subset R^{3}$, and\ $\Omega$\ is bounded,\ $\partial\Omega$\ is smooth or piecewise smooth, \[u=(u_{1}(x,\ y,\ z,\ t),\cdots,u_{m}(x,\ y,\ z,\ t))^{T}\in C^{1}(\Omega\times (0,\ T)),\] the initial conditions and boundary conditions are\[u|_{t=0}\in L^{2}(\Omega),\ u|_{\partial\Omega\times(0,\ T)}\in L^{2}(\partial\Omega\times(0,\ T)),\] $f_{j},\ 1\leq j\leq m,$\ are continuous functions? Of course it isn't. We should solve Eq.(4.1) for\ $m$\ components of\ $u_{t},\ u,\ u_{x},\ u_{y},\ u_{z}$.\\Let's consider an example of the second order equation at first,\[u_{t}=f(u,\ u_{x},\ u_{y},\ u_{z},\ u_{xx},\ u_{xy},\ u_{xz},\ u_{yy},\ u_{yz},\ u_{zz},\ x,\ y,\ z,\ t). \] We should transform it into the first order equations as follows,\begin{eqnarray}&&u_{t}=f(u,\ u_{x},\ u_{y},\ u_{z},\ v_{x},\ v_{y},\ v_{z},\ w_{y},\ w_{z},\ r_{z},\ x,\ y,\ z,\ t),\\&&v=u_{x},\\&&w=u_{y},\\&&r=u_{z}.\end{eqnarray}We notice that it is not always the derivatives in the left hand sides. It should be resolved with any derivatives in the left hand side, including the components of\ $u$. Hence the general first order partial differential equations should be like the following,\[v_{j}=f_{j}(v_{m+1},\ v_{m+2},\cdots,\ v_{5m},\ x,\ y,\ z,\ t),\ 1\leq j\leq m,\] where\ $t\in [0,\ T],\ (x,\ y,\ z)^{T}\in \Omega\subset R^{3}$, and\ $\Omega$\ is bounded,\ $\partial\Omega$\ is smooth or piecewise smooth, \[u=(u_{1}(x,\ y,\ z,\ t),\cdots,u_{m}(x,\ y,\ z,\ t))^{T}\in C^{1}(\Omega\times (0,\ T)).\] The initial conditions and boundary conditions are\[u|_{t=0}\in L^{2}(\Omega),\ u|_{\partial\Omega\times(0,\ T)}\in L^{2}(\partial\Omega\times(0,\ T)).\] $f_{j},\ 1\leq j\leq m,$\ are continuous functions,\ $v_{1},\ v_{2},\cdots,\ v_{5m}$\ is a permutation of the components of\ $u_{t},\ u,\ u_{x},\ u_{y},\ u_{z}.$\\We assume there are\ $r$\ components of\ $u$\ in the beginning of the right hand sides of the Eq.(4.9)\[v_{m+1},\ v_{m+2},\cdots,\ v_{m+r},\ 0\leq r\leq m.\] And there are\ $4m-r$\ derivatives of\ $u$\ in the right hand sides of the Eq.(4.9). This point is very important in the following.\\We transform Eq.(4.9) into the
linear equations on unknown functions as follows,\[v_{j}-\sum_{k=1}^{4m}c_{jk}v_{m+k}-s_{j}=0,\ 1\leq j\leq m,\]where\ $c_{jk},\ 1\leq j\leq m,\ 1\leq k\leq 4m,$\ are all real constants to be determined,\[s_{j}=f_{j}(v_{m+1},\ v_{m+2},\cdots,\ v_{5m},\ x,\ y,\ z,\ t)-\sum_{k=1}^{4m}c_{jk}v_{m+k},\ 1\leq j\leq m.\]
Let's introduce\ $X=(V_{1}^{T},\ V_{2}^{T},\ V_{3}^{T},\ V_{4}^{T},\ V_{5}^{T},\ S^{T})^{T}$, where\[V_{i}=(v_{m(i-1)+1},\ v_{m(i-1)+2},\cdots,\ v_{mi})^{T},\ 1\leq i\leq5,\ S=(s_{1},\ s_{2},\cdots,\ s_{m})^{T}.\]Then Eq.(4.13) is equivalent to\[\alpha^{T}X=0,\] where\[ \alpha^{T}=(E,\ -C_{1},\ -C_{2},\ -C_{3},\ -C_{4},\ -E),\ C_{i}=(c_{j,(m(i-1)+k)})_{m\times m},\ 1\leq i\leq 4.\]
 After we solve Eq.(4.16) by linear algebra, we obtain\ $\exists$\ independent variable vector\ $Z$, such that\[X=\left(
                                                                                                                                     \begin{array}{c}
                                                                                                                                       \beta^{T} \\
                                                                                                                                       E_{5m} \\
                                                                                                                                     \end{array}
                                                                                                                                   \right)Z,\ \mbox{where}\ \beta^{T}=(C_{1},\ C_{2},\ C_{3},\ C_{4},\ E).\]
                                                                                                                                   We should discuss the independent variable vector\ $Z$\ as follows,
                                                                                                                                   \begin{eqnarray}&&u_{t}=A_{00}Z,\ u=A_{0}Z,\ u_{x}=A_{01}Z,\ u_{y}=A_{02}Z,\ u_{z}=A_{03}Z,\ S=E_{5}^{T}Z,\end{eqnarray} where\ $E_{j}=(e_{m(j-1)+1},\ e_{m(j-1)+2},\cdots,\ e_{mj}),\ 1\leq j\leq 5,$\
 and the rows of\ $A_{00},\ A_{0},\ A_{01},\ A_{02},\ A_{03}$\ are a permutation of the rows of\ $\beta^{T},\ E_{1}^{T},\ E_{2}^{T},\ E_{3}^{T},\ E_{4}^{T}$.\\
 And we obtain Eq.(4.9) is equivalent to the following system respect to\ $Z$,\begin{eqnarray}&&\cfrac{\partial A_{0}Z}{\partial t}=A_{00}Z,\\&& \cfrac{\partial A_{0}Z}{\partial x}=A_{01}Z,\\&& \cfrac{\partial A_{0}Z}{\partial y}=A_{02}Z,\\&& \cfrac{\partial A_{0}Z}{\partial z}=A_{03}Z,\\&& f(E_{1}^{T}Z,E_{2}^{T}Z,E_{3}^{T}Z,E_{4}^{T}Z,x, y, z, t)-C_{1}E_{1}^{T}Z-C_{2}E_{2}^{T}Z-C_{3}E_{3}^{T}Z-C_{4}E_{4}^{T}Z=E_{5}^{T}Z.\end{eqnarray}
 We will also solve Eq.(4.20) to Eq.(4.23) by the Fourier transform as follows. Then Eq.(4.24) is the goal. \\
 \begin{eqnarray}\int_{0}^{T}dt\int_{\Omega}e^{-it\xi_{0}-ix\xi_{1}-iy\xi_{2}-iz\xi_{3}}\cfrac{\partial A_{0}Z}{\partial t}dxdydz=\int_{0}^{T}dt\int_{\Omega}e^{-it\xi_{0}-ix\xi_{1}-iy\xi_{2}-iz\xi_{3}}A_{00}Zdxdydz,&&\\ \int_{0}^{T}dt\int_{\Omega}e^{-it\xi_{0}-ix\xi_{1}-iy\xi_{2}-iz\xi_{3}} \cfrac{\partial A_{0}Z}{\partial x}dxdydz=\int_{0}^{T}dt\int_{\Omega}e^{-it\xi_{0}-ix\xi_{1}-iy\xi_{2}-iz\xi_{3}}A_{01}Zdxdydz,&&\\ \int_{0}^{T}dt\int_{\Omega}e^{-it\xi_{0}-ix\xi_{1}-iy\xi_{2}-iz\xi_{3}}\cfrac{\partial A_{0}Z}{\partial y}dxdydz=\int_{0}^{T}dt\int_{\Omega}e^{-it\xi_{0}-ix\xi_{1}-iy\xi_{2}-iz\xi_{3}}A_{02}Zdxdydz,&&\\ \int_{0}^{T}dt\int_{\Omega}e^{-it\xi_{0}-ix\xi_{1}-iy\xi_{2}-iz\xi_{3}} \cfrac{\partial A_{0}Z}{\partial z}dxdydz=\int_{0}^{T}dt\int_{\Omega}e^{-it\xi_{0}-ix\xi_{1}-iy\xi_{2}-iz\xi_{3}}A_{03}Zdxdydz.&&\end{eqnarray}\\
By using the definition 2.1, we obtain\begin{eqnarray}
FI(\cfrac{\partial A_{0}Z}{\partial t})&=&\int_{0}^{T}dt\int_{\Omega}e^{-it\xi_{0}-ix\xi_{1}-iy\xi_{2}-iz\xi_{3}}(\cfrac{\partial A_{0}Z}{\partial t})dxdydz
\\&=&\int_{\Omega}(A_{0}Z)e^{-it\xi_{0}}|_{t=0}^{t=T}e^{-ix\xi_{1}-iy\xi_{2}-iz\xi_{3}}dxdydz
+\\&&i\xi_{0}\int_{0}^{T}dt\int_{\Omega}e^{-it\xi_{0}-ix\xi_{1}-iy\xi_{2}-iz\xi_{3}}(A_{0}Z)dxdydz
\\&=&g_{0}+ i\xi_{0} FI(A_{0}Z),\end{eqnarray}\begin{eqnarray}
FI(\cfrac{\partial A_{0}Z}{\partial x})&=&\int_{0}^{T}dt\int_{\Omega}e^{-it\xi_{0}-ix\xi_{1}-iy\xi_{2}-iz\xi_{3}}(\cfrac{\partial A_{0}Z}{\partial x})dxdydz
\\&=&\int_{0}^{T}\int_{\partial\Omega}(A_{0}Z)n_{1}e^{-it\xi_{0}-ix\xi_{1}-iy\xi_{2}-iz\xi_{3}}dxdydz
+\\&&i\xi_{1}\int_{0}^{T}dt\int_{\Omega}e^{-it\xi_{0}-ix\xi_{1}-iy\xi_{2}-iz\xi_{3}}(A_{0}Z)dxdydz
\\&=&g_{1}+ i\xi_{1} FI(A_{0}Z),\end{eqnarray} \begin{eqnarray}
FI(\cfrac{\partial A_{0}Z}{\partial y})&=&\int_{0}^{T}dt\int_{\Omega}e^{-it\xi_{0}-ix\xi_{1}-iy\xi_{2}-iz\xi_{3}}(\cfrac{\partial A_{0}Z}{\partial y})dxdydz
\\&=&\int_{0}^{T}\int_{\partial\Omega}(A_{0}Z)n_{2}e^{-it\xi_{0}-ix\xi_{1}-iy\xi_{2}-iz\xi_{3}}dxdydz
+\\&&i\xi_{2}\int_{0}^{T}dt\int_{\Omega}e^{-it\xi_{0}-ix\xi_{1}-iy\xi_{2}-iz\xi_{3}}(A_{0}Z)dxdydz
\\&=&g_{2}+ i\xi_{2} FI(A_{0}Z),\end{eqnarray} \begin{eqnarray}
FI(\cfrac{\partial A_{0}Z}{\partial z})&=&\int_{0}^{T}dt\int_{\Omega}e^{-it\xi_{0}-ix\xi_{1}-iy\xi_{2}-iz\xi_{3}}(\cfrac{\partial A_{0}Z}{\partial z})dxdydz
\\&=&\int_{0}^{T}\int_{\partial\Omega}(A_{0}Z)n_{3}e^{-it\xi_{0}-ix\xi_{1}-iy\xi_{2}-iz\xi_{3}}dxdydz
+\\&&i\xi_{3}\int_{0}^{T}dt\int_{\Omega}e^{-it\xi_{0}-ix\xi_{1}-iy\xi_{2}-iz\xi_{3}}(A_{0}Z)dxdydz
\\&=&g_{3}+ i\xi_{3} FI(A_{0}Z),\end{eqnarray}where\begin{eqnarray}g_{0}&=&\int_{\Omega}(A_{2}e^{-iT\xi_{0}}-A_{1})e^{-ix\xi_{1}-iy\xi_{2}-iz\xi_{3}}dxdydz,\\
g_{1}&=&\int_{0}^{T}dt\int_{\partial \Omega}A_{3}n_{1}e^{-it\xi_{0}-ix\xi_{1}-iy\xi_{2}-iz\xi_{3}}dS,\\
g_{2}&=&\int_{0}^{T}dt\int_{\partial \Omega}A_{3}n_{2}e^{-it\xi_{0}-ix\xi_{1}-iy\xi_{2}-iz\xi_{3}}dS,\\
g_{3}&=&\int_{0}^{T}dt\int_{\partial \Omega}A_{3}n_{3}e^{-it\xi_{0}-ix\xi_{1}-iy\xi_{2}-iz\xi_{3}}dS.\\
A_{1}&=&u|_{t=0},\ A_{2}=u|_{t=T},\ A_{3}=u|_{\partial\Omega\times(0,\ T)},
\end{eqnarray}\ $n_{k}$\ is the\ $k$th component of
the normal vector to\ $\partial \Omega,\ k=1,\ 2,\ 3$. We only need\[A_{1},\ A_{2}\in L^{2}(\Omega),\ A_{3}\in L^{2}(\partial\Omega\times(0,\ T)).\]
Now we transformed the equations Eq.(4.20) to Eq.(4.23) into the following.\[\label
{A-Problem}BFI(Z)=\beta_{1},\]where\[B=\left(
                                         \begin{array}{c}
                                        i\xi_{0}A_{0}-A_{00} \\
                                         i\xi_{1}A_{0}-A_{01} \\
                                         i\xi_{2}A_{0}-A_{02} \\
                                         i\xi_{3}A_{0}-A_{03} \\
                                         \end{array}
                                       \right)_{4m\times5m}
=(B_{1},\ -B_{2}),\ \beta_{1}=(-g_{0}^{T},\ -g_{1}^{T},\ -g_{2}^{T},\ -g_{3}^{T})^{T},\] $B_{1}$\ is the first\ $4m$\ columns of\ $B$,\ $-B_{2}$\ is the last\ $m$\ columns of\ $B$.\\Because the rows of\ $A_{00},\ A_{0},\ A_{01},\ A_{02},\ A_{03}$\ are the permutation of the rows of\ $\beta^{T},\ E_{4m}$, we obtain the rows of\ $B_{2}$\ are the permutation of the rows of\ $E$\ and\ $0$. And we assumed there are\ $r$\ components of\ $u$\ in the beginning of the right hand sides of the Eq.(4.9)\[v_{m+1},\ v_{m+2},\cdots,\ v_{m+r},\ 0\leq r\leq m.\] There are\ $4m-r$\ derivatives of\ $u$\ in the right hand sides of the Eq.(4.9). Hence there exists\ $P$\ is a permutation matrix, such that\[P\left(
                                  \begin{array}{c}
                                    A_{00}\\
                                    A_{01} \\
                                    A_{02}\\
                                    A_{03} \\
                                  \end{array}
                                \right)_{B_{1}}=\left(
                                          \begin{array}{cc}
                                            C_{r1}, &C_{r2} \\
                                            0, & E_{4m-r}\\
                                          \end{array}
                                        \right),\ \mbox{where}\ \left(
                                  \begin{array}{c}
                                    A_{00}\\
                                    A_{01} \\
                                    A_{02}\\
                                    A_{03} \\
                                  \end{array}
                                \right)_{B_{1}}\ \mbox{is the first\ $4m$\ columns of}\ \left(
                                  \begin{array}{c}
                                    A_{00}\\
                                    A_{01} \\
                                    A_{02}\\
                                    A_{03} \\
                                  \end{array}
                                \right).\]
Therefore, we obtain\[PB_{1}=\left(
                               \begin{array}{cc}
                                 A_{r1}-C_{r1},&  A_{r2}-C_{r2} \\
                                A_{r3}, & A_{r4}-E_{4m-r} \\
                               \end{array}
                             \right),\]where\ $A_{r1},\ A_{r2},\ A_{r3},\ A_{r4}$\ are related with\ $i\xi_{0},\ i\xi_{1},\ i\xi_{2},\ i\xi_{3}$, the elements in them are not constants except\ $0$. If\ $i\xi_{0}=i\xi_{1}=i\xi_{2}=i\xi_{3}=0$, then\ $ A_{r4}=0$. Hence\ $det(A_{r4}-E_{4m-r})$\ is not always\ $0$. And\ $A_{r4}-E_{4m-r}$\ is convertible. We obtain
                             \[det(PB_{1})=det(A_{r1}-C_{r1}-(A_{r2}-C_{r2})(A_{r4}-E_{4m-r})^{-1}A_{r3})det(A_{r4}-E_{4m-r}).\]
                             Now we determine the parameters. We choose\ $C_{1},\ C_{2},\ C_{3},\ C_{4}$, such that\ $det(PB_{1})\neq 0$. Hence\ $det(B_{1})\neq 0,\ F^{-1}(B_{1}^{-1})$\ and\ $F^{-1}(B_{1}^{-1}B_{2})$\ exists.\\
Now we assume
\[ Z=\left(
       \begin{array}{c}
         Z_{1} \\
         Z_{2} \\
       \end{array}
     \right),\ \mbox{where}\ Z_{1}\ \mbox{is the first\ $4m$\ componenets of}\ Z,\] then we obtain\[B_{1}FI(Z_{1})=\beta_{1}+B_{2}FI(Z_{2}).\]
                                                                                                                                   If we assume\ $C_{0}=\{\xi^{\prime}|det(B_{1})=0\}$, where\ $\xi^{\prime}=(\xi_{0},\ \xi_{1},\ \xi_{2},\ \xi_{3})^{T}$, then the measure of\ $C_{0}$\ is\ $0$. And we obtain
 \[FI(Z_{1})(1-I_{C_{0}}(\xi^{\prime}))=B_{1}^{-1}\beta_{1}(1-I_{C_{0}}(\xi^{\prime}))+B_{1}^{-1}B_{2}FI(Z_{2})(1-I_{C_{0}}(\xi^{\prime})).\]
From two lemmas in second section, we also obtain\[F^{-1}[FI(Z_{1})(1-I_{C_{0}}(\xi^{\prime}))]=Z_{1}I_{\Omega\times(0,\ T)},\ F^{-1}[FI(Z_{2})(1-I_{C_{0}}(\xi^{\prime}))]=Z_{2}I_{\Omega\times(0,\ T)}.\]
If we assume\[w_{1}(x,\ y,\ z,\ t)=F^{-1}[B_{1}^{-1}\beta_{1}(1-I_{C_{0}}(\xi^{\prime}))],\ w_{2}(x,\ y,\ z,\ t)=F^{-1}(B_{1}^{-1}B_{2}),\]
then we obtain\[Z_{1}I_{\Omega\times(0,\ T)}=w_{1}(x,\ y,\ z,\ t)+w_{2}(x,\ y,\ z,\ t)\ast(Z_{2}I_{\Omega\times(0,\ T)}),\]
where\begin{eqnarray}Z_{1}&=&(E_{1}^{T}Z,\ E_{2}^{T}Z,\ E_{3}^{T}Z,\ E_{4}^{T}Z)^{T}=(V_{1}^{T},\ V_{2}^{T},\ V_{3}^{T},\ V_{4}^{T})^{T},\\Z_{2}&=&E_{5}^{T}Z=f(E_{1}^{T}Z,\ E_{2}^{T}Z,\ E_{3}^{T}Z,\ E_{4}^{T}Z,\ x,\ y,\ z,\ t)-\\&& C_{1}E_{1}^{T}Z-C_{2}E_{2}^{T}Z-C_{3}E_{3}^{T}Z-C_{4}E_{4}^{T}Z=S.\end{eqnarray}It is obvious\ $\exists\ \psi$, such that\ $Z_{2}=\psi(Z_{1})$. Therefore, we obtain
\[ Z_{1}I_{\Omega\times(0,\ T)}=w_{1}(x,\ y,\ z,\ t)+w_{2}(x,\ y,\ z,\ t)\ast(\psi(Z_{1})I_{\Omega\times(0,\ T)}).\]
If\ $Z_{1}$\ satisfied Eq.(4.66), then we let\ $Z_{2}=\psi(Z_{1})$. We obtain\ $ B FI(Z)=\beta_{1},\ \alpha^{T}X=0,$\ on\ $\Omega\times(0,\ T)$. Therefore,\ $A_{0}Z$\ is the solution of Eq.(4.9)\ on\ $\Omega\times(0,\ T)$. Hence we arrive at \begin{theorem} \label{Theorem2-1}\ $ w_1,\ w_2,\ \psi,$\ as we described, then Eq.(4.9) is equivalent to Eq.(4.66). \end{theorem}
At this time, we will also show that we should't know\ $A_{2}=u|_{t=T}$. We can get it by\begin{eqnarray}&&\int_{-\infty}^{+\infty}e^{-iT\xi_{0}} e^{it\xi_{0}}d\xi_{0}=\delta(t-T),\\&& F^{-1}(B_{1}^{-1}e^{-iT\xi_{0}})=(\varphi(x,\ y,\ z,\ t)I_{\{t\geq 0\}})\ast\delta(t-T)=\varphi(x,\ y,\ z,\ t-T)I_{\{t\geq T\}}.\end{eqnarray}
We only need to prove\[F^{-1}(B_{1}^{-1})=F^{-1}(B_{1}^{-1})I_{\{t\geq 0\}}. \]From\ $B_{1}^{-1}=B_{1}^{\ast}/det(B_{1})$, and\ $B_{1}^{\ast}$\ is a polynomial matrix, we only need to prove\[F^{-1}((det(B_{1}))^{-1})=F^{-1}((det(B_{1}))^{-1})I_{\{t\geq 0\}}. \]
We can transform\ $B_{1}$\ into the following by the primary row transformations on the rows from\ $(m+1)th$\ to\ $(4m)th$\ and some transpositions of the columns,\[ \left(
            \begin{array}{cccc}
              i\xi_{0}A_{11}-B_{11}, & i\xi_{0}A_{12}-B_{12}, & i\xi_{0}A_{13}-B_{13}, &i\xi_{0}A_{14}-B_{14}, \\
              C_{21} & E & 0 & 0 \\
              C_{31} &0 & E & 0 \\
              C_{41} & 0 & 0 & E \\
            \end{array}
          \right),\]where\ $C_{21},\ C_{31},\ C_{41}$\ are not related to\ $i\xi_{0},\ A_{1j},\ B_{1j},\ 1\leq j\leq 4,$\ are all constant matrices. Hence we obtain
          \begin{eqnarray}&&det(B_{1})=det(i\xi_{0}A_{05}-B_{05})\phi(i\xi_{1},\ i\xi_{2},\ i\xi_{3}),\\&& \mbox{where}\ A_{05}=A_{11}-\sum_{j=2}^{4}A_{1j}C_{j1},\ B_{05}=B_{11}-\sum_{j=2}^{4}B_{1j}C_{j1}.\end{eqnarray}$A_{05}$\ and\ $B_{05}$\ are not related to\ $i\xi_{0}$, neither. Then we choose the parameter matrices\ $C_{1},\ C_{2},\ C_{3},\ C_{4}$\ which satisfy\ $A_{05}$\ is convertible and\ $Re(\lambda)<0$, where\ $\lambda$\ is the characteristic value of matrix\ $A_{05}^{-1}B_{05}$. This is available because there are\ $4m^{2}$\ variables in the parameter matrices\ $C_{1},\ C_{2},\ C_{3},\ C_{4}$\ and\ $A_{05},\ B_{05}$\ are\ $m\times m$. From the remark 3.1, we obtain\begin{eqnarray}&&det(B_{1})=det(A_{05})det(i\xi_{0}E-A_{05}^{-1}B_{05})\phi(i\xi_{1},\ i\xi_{2},\ i\xi_{3}),\\&&\int_{-\infty}^{+\infty}(det(i\xi_{0}E-A_{05}^{-1}B_{05}))^{-1}e^{it\xi_{0}}d\xi_{0}=\int_{-\infty}^{+\infty}(det(i\xi_{0}E-A_{05}^{-1}B_{05}))^{-1}e^{it\xi_{0}}d\xi_{0}I_{\{t\geq 0\}},\\&& \int_{-\infty}^{+\infty}(det(B_{1}))^{-1}e^{it\xi_{0}}d\xi_{0}=\int_{-\infty}^{+\infty}(det(B_{1}))^{-1}e^{it\xi_{0}}d\xi_{0}I_{\{t\geq 0\}}.\end{eqnarray} Hence (4.69) stands. We also take\ $g_{0}$\ as\[g_{0}=\int_{\Omega}(-A_{1})e^{-ix\xi_{1}-iy\xi_{2}-iz\xi_{3}}dxdydz.\]Now we have transformed the first order partial differential equations resolved with any derivatives into the equivalent integral equations as Eq.(4.66).
         \section{Classical solution and generalized solution}\setcounter{equation}{0}
We have transformed the first order partial differential equations resolved with any derivatives into the equivalent integral equations as follows,
\[ Z_{1}I_{\Omega\times(0,\ T)}=w_{1}(x,\ y,\ z,\ t)+w_{2}(x,\ y,\ z,\ t)\ast(\psi(Z_{1})I_{\Omega\times(0,\ T)}),\]
where\[ Z_{1}=(Z_{1j})_{p\times1},\ w_{1}=(w_{1j})_{p\times1},\ w_{2}=(w_{2}^{i,j})_{p\times q},\ \psi(Z_{1})=(\psi_{j}(Z_{1}))_{q\times1},\] $p,\ q$ are natural numbers. We notice that\[w_{1}(x,\ y,\ z,\ t)\neq w_{1}(x,\ y,\ z,\ t)I_{\Omega\times(0,\ T)},\ w_{2}(x,\ y,\ z,\ t)\neq w_{2}(x,\ y,\ z,\ t)I_{\Omega\times(0,\ T)}.\] Then (5.1) is equivalent to\begin{eqnarray}&&Z_{1}I_{\Omega\times(0,\ T)}=w_{1}(x,\ y,\ z,\ t)+w_{2}(x,\ y,\ z,\ t)\ast(\psi(Z_{1})I_{\Omega\times(0,\ T)}),\ \forall (x,\ y,\ z,\ t)\in \Omega_{0},\\&&0=w_{1}(x,\ y,\ z,\ t)+w_{2}(x,\ y,\ z,\ t)\ast(\psi(Z_{1})I_{\Omega\times(0,\ T)}),\ \mbox{otherwise},\end{eqnarray}where\ $\Omega_{0}=\Omega\times(0,\ T)$.\\If there exists\ $Z_{1}I_{\Omega\times(0,\ T)}\in C(\Omega\times(0,\ T))$\ satisfies Eq.(5.4) and Eq.(5.5) both, then we say the classical solution of the first order partial differential equations resolved with any derivatives exists. If there exists\ $Z_{1}I_{\Omega\times(0,\ T)}\in L^{2}(\Omega\times(0,\ T))$\ only satisfies Eq.(5.4), then we say the generalized solution of the first order partial differential equations resolved with any derivatives exists. Maybe Eq.(5.5) can explain why sometimes the classical solution doesn't exist. The generalized solution is always locally exist and unique.
 \section{Acknowledgements}
 I give my best thanks to my supervisor Prof. Mark Edelman, for his guidance when I am a Scholar Visitor in Yeshiva University. I sincerely thank Prof. Caisheng Chen in Hohai University, Prof. Junxiang Xu in Southeast University and Prof. Zuodong Yang in Nanjing Normal University for their recommendation and other helps.\\The financial support of Chinese ministry of education is gratefully acknowledged. \\

\end{document}